\documentclass[12pt]{article}
\font\smallit=cmti10

\usepackage{amsmath}
\usepackage{amssymb}
\usepackage{amsthm}
\usepackage[dvips]{graphicx}
\usepackage{color}
\usepackage{epsfig}
\usepackage[mathscr]{eucal}
\newtheorem{theorem}{Theorem}

\newtheorem{lemma}[theorem]{Lemma}
\newtheorem{cor}[theorem]{Corollary}

\allowdisplaybreaks

\begin{document} 

\begin{center}
{\bf CONTINUED FRACTIONS WITH PARTIAL QUOTIENTS BOUNDED IN AVERAGE}
\vskip 20pt
{\bf Joshua N. Cooper}\\
{\smallit ETH-Z\"urich, Institute for Theoretical Computer Science}\\
{\smallit Universit\"atstr.\@ 6, Z\"urich, Switzerland CH-8092}\\
{\tt cooper@cims.nyu.edu}\\
\vskip 10pt
\end{center}
\vskip 30pt 


\begin{abstract}
We ask, for which $n$ does there exists a $k$, $1 \leq k < n$ and $(k,n)=1$, so that $k/n$ has a continued fraction whose partial quotients are bounded in average by a constant $B$?  This question is intimately connected with several other well-known problems, and we provide a lower bound in the case of $B=2$.  The proof, which is completely elementary, involves a simple ``shifting'' argument, the Catalan numbers, and the solution to a linear recurrence.
\end{abstract}

\thispagestyle{empty} 
\baselineskip=15pt 
\vskip 30pt 

\section{Introduction}

An important question in the theory of quasirandomness, uniform distribution of points, and diophantine approximation is the following: For which $n \in \mathbb{Z}$ is it true that there exists an integer $k$, $1 \leq k < n$ and $(k,n)=1$, so that $k/n$ has a continued fraction whose partial quotients are bounded in average by a constant $B$?  That is, if we write $k/n = [0;a_1,a_2,\ldots,a_m]$, we wish to find $k$ so that
$$
t^{-1} \sum_{i=1}^t a_i \leq B
$$
for all $t$ with $1 \leq t \leq m$.  Denote by $\mathcal{F}(B)$ the set of all $n$ for which such a $k$ exists.  These sets are discussed at length in \cite{C2} and the related matter of partial quotients bounded {\it uniformly} by a constant appears as an integral part of \cite{N1}.  This latter question is closely connected with Zaremba's Conjecture (\cite{Z}), which states that such a $k$ exists for all $n > 1$ if we take $B=5$.

Define the {\it continuant} $K(a_1,a_2,\ldots,a_m)$ to be the denominator of the continued fraction $k/n = [0;a_1,a_2,\ldots,a_m]$. In \cite{Cu1}, it is proven that, if $S_n(B)$ is the number of sequences $\mathbf{a}=(a_1,\ldots,a_m)$ bounded uniformly by $B$ with $K(\mathbf{a}) \leq n$ and $H(B)$ is the Hausdorff dimension of the set of continued fractions with partial quotients bounded uniformly by $B$, then
$$
\lim_{n \rightarrow \infty} \frac{\log(S_n(B))}{\log n} = 2 H(B).
$$
Then, in \cite{He1}, $H(2)$ is calculated with a great deal of accuracy: $H(2) \approx 0.53128$.  Therefore, $S_n(2)$, and thus the number of $p/q$ with $q \leq n$ whose partial quotients are bounded by $2$, is $n^{1.0625\ldots + o(1)}$.  (This improves the previous best known lower bound, $n^{\approx 1.017}$ computed in \cite{Cu1}, slightly.)

Define $\bar{S}_n(B)$ to be the number of sequences $\mathbf{a} = (a_1,\ldots,a_m)$ with partial quotients bounded {\it in average} by $B$ so that $K(\mathbf{a}) \leq n$.  Clearly, $\bar{S}_n(B) \geq S_n(B)$, so $\bar{S}_n(2) \gg n^{1.0625}$.  In the next section, we prove something much stronger, however -- an exponent of $\approx 1.5728394$ -- thus providing a lower bound in the first nontrivial case.  Section \ref{F2} discusses the implications for the density of $\mathcal{F}(2)$ and a few open problems.

\vskip 30pt

\section{The Proof}

\begin{theorem} \label{newbnd} For any $\epsilon > 0$, $\bar{S}_n(2) \gg n^{2\log 2/\log(1+\sqrt{2})-\epsilon}$.
\end{theorem}
\begin{proof} The proof consists of two parts: computing the number of positive sequences of length $m$ bounded in average by $2$, and then computing the smallest possible $m$ so that $K(a_1,\ldots,a_m) > n$ and the $a_i$ are bounded in average by $2$.

First, we wish to know how many sequences $(a_1,\ldots,a_m)$ there are with $a_j \geq 1$ for each $j \in [m]$ and $\sum_{j=1}^r a_j \leq 2r$ for each $r \in [m]$.  Call this number $T(m)$.  By writing $b_j = a_j - 1$, we could equivalently ask for sequences $(b_1,\ldots,b_m)$ with $b_j \geq 0$ for each $j \in [m]$ and $\sum_{j=1}^r b_j \leq r$ for each $r \in [m]$.  This is precisely the number of lattice paths from $(0,0)$ to $(m,m)$ which do not cross the line $y=x$, and so $T(m)$ is the $m^\text{th}$ Catalan number, or $(m+1)^{-1} \binom{2m}{m} = 4^{m(1-o(1))}$.

In the following lemmas, we show that $K(a_1,\ldots,a_m) \leq n$ if $m \leq \log n (1-o(1))/\log(1+\sqrt{2})$.  Therefore, setting $m$ as large as possible, we have at least
$$
4^{\log n (1-o(1)) / \log(1+\sqrt{2})} = n^{2\log 2/\log(1+\sqrt{2})-o(1)}
$$
sequences with partial quotients bounded in average by $2$ and continuant $\leq n$ .
\end{proof}

We must show that the size of a continuant with partial quotients bounded in average by $B$ is at most the largest size of a continuant with partial quotients bounded by $B$.

\begin{lemma} If the sequence $(a_1,\ldots,a_m)$ of positive integers is bounded in average by $B>1$, then $K(a_1,\ldots,a_m) \leq K(\underbrace{B,\ldots,B}_m)$.
\end{lemma}
\begin{proof} We prove the Lemma by a ``shifting'' argument.  That is, we perform induction on the size of the entry $a_j$ such that $a_j>B$ and $j$ is as small as possible.  If $\mathbf{a}=(a_1,\ldots,a_m)$ contains no $a_t > B$, we are done, because increasing the partial quotients can only increase the continuant.  If there is some $a_t > B$, let $t \geq 2$ be the smallest such index.  We consider two cases: (i) $a_t \geq B+2$ or $a_{t-1} < B$, and (ii) $a_t = B+1$, $a_k = B$ for $s \leq k \leq t-1$ for some $2 \leq s \leq t-1$, and $a_{s-1} < B$.  (Clearly, $\mathbf{a} \neq (B,B,\ldots,B,B+1,a_{t+1},\ldots,a_m)$, since this sequence is not bounded in average by $B$.  Therefore we may assume $s\geq 2$.)

\noindent Case (i):

Let $\mathbf{b}=(b_1,\ldots,b_m)$ = $(a_1,\ldots,a_{t-1}+1,a_t-1,\ldots,a_m)$.  We show that $K(\mathbf{b}) > K(\mathbf{a})$.  First, note that
$$
\sum_{j=1}^r b_j = \left \{ \begin{array}{ll} \sum_{j=1}^r a_j & \text{ if } r \neq t-1 \\ 1+\sum_{j=1}^{t-1} a_j & \text{ if } r=t-1. \end{array} \right .
$$
Since $a_t \geq B+1$, $\sum_{j=1}^{t-1} a_j \leq tB-B-1$, so $1+\sum_{j=1}^{t-1} a_j \leq (t-1)B$, and $\mathbf{b}$ is bounded in average by $B$.  Second, note that it suffices to consider the case of $t=m$, since, if $K(b_1,\ldots,b_{j}) > K(a_1,\ldots,a_{j})$ for $1 \leq j \leq t$, then $K(\mathbf{b}) > K(\mathbf{a})$.  (That is, $K(\cdot)$ is monotone increasing.)

Let $q_j = K(a_1,\ldots,a_{j})$ and $q^\prime_j = K(b_1,\ldots,b_{j})$.  (We use the convention that $q_j = 0$ when $j < 0$ and $q_0 = 1$.)  Clearly, $q_j = q^\prime_j$ if $j < t-1$.  When $j=t-1$, we have $q^\prime_{t-1} > q_{t-1}$ by monotonicity.  When $j=t$,
$$
q_t = a_t q_{t-1} + q_{t-2} = a_t (a_{t-1} q_{t-2} + q_{t-3}) + q_{t-2} = (a_t a_{t-1} + 1) q_{t-2} + a_t q_{t-3},
$$
and
\begin{align*}
q^\prime_t &= (b_t b_{t-1} + 1) q^\prime_{t-2} + b_t q^\prime_{t-3} \\
&= ((a_t-1)(a_{t-1}+1)+1) q_{t-2} + (a_t-1) q_{t-3} \\
&= q_t + q_{t-2} (a_t - a_{t-1} - 1) - q_{t-3}.
\end{align*}
Since $a_t \geq a_{t-1} + 2$ and $q_{t-2} > q_{t-3}$, we have
$$
q^\prime_t \geq q_t + q_{t-2} - q_{t-3} > q_t.
$$

\noindent Case (ii).

Now, assume that $a_t = B+1$, $a_k = B$ for $s \leq k \leq t-1$ for some $2 \leq s \leq t-1$, and $a_{s-1} < B$.  Then define $\mathbf{b} = (b_1,\ldots,b_m)$ by letting $b_j = a_j$ if $j \neq s-1$ and $j \neq t$; $b_{s-1} = a_{s-1} + 1$; and $b_t = a_t - 1$.  Again, we may assume that $t=m$.  Then
$$
\sum_{j=1}^r b_j = \left \{ \begin{array}{ll} \sum_{j=1}^r a_j & \text{ if } r = t \text{ or } r < s-1 \\ 1+\sum_{j=1}^{r} a_j & \text{ if } s-1 \leq r \leq t-1. \end{array} \right .
$$
For any $r$ such that $s-1 \leq r \leq t-1$,
$$
\sum_{j=1}^r a_j = \sum_{j=1}^{t} a_j - \sum_{j=r+1}^t a_j \leq Bt - (B(t-r-1)+(B+1)) \leq Br -1.
$$
Therefore, $\sum_{j=1}^r b_j \leq Br$ for all $r \in [t]$, and we may conclude that $\mathbf{b}$ is bounded in average by $B$.

Define $F_k$ as follows: $F_0=0$, $F_1=1$, and, for $k>1$, $F_k=BF_{k-1}+F_{k-2}$.  Then it is easy to see by induction that
$$
K(\underbrace{B,\ldots,B}_k,x)=F_{k+1}x + F_k.
$$
Also,
\begin{equation} \label{recursion}
K(y,c_1,\ldots,c_r) = y K(c_1,\ldots,c_r) + K(c_2,\ldots,c_r).
\end{equation}
Taking $k = t-s$, we deduce
$$
K(a_{s-1},\ldots,a_t) = a_{s-1} ((B+1) F_{k+1} + F_k) + (B+1) F_k + F_{k-1},
$$
and
\begin{align*}
K(b_{s-1},\ldots,b_t) &= (a_{s-1}+1) (BF_{k+1} + F_k) + BF_k + F_{k-1} \\
&= K(a_{s-1},\ldots,a_t) + (B - a_{s-1}) F_{k+1} \\
& \geq K(a_{s-1},\ldots,a_t) + F_{k+1}.
\end{align*}
If $s=2$, we are done.  Otherwise, we use that
\begin{align*}
K(b_{s-2},\ldots,b_t) &= a_{s-2} K(b_{s-1},\ldots,b_t) + K(b_s,\ldots,b_t) \\
&\geq a_{s-2} K(a_{s-1},\ldots,a_t) + F_{k+1} + b_t K(b_s,\ldots,b_{t-1}) + K(b_s,\ldots,b_{t-2})\\
&= a_{s-2} K(a_{s-1},\ldots,a_t) + F_{k+1} + K(a_s,\ldots,a_t) - K(a_s,\ldots,a_{t-1}) \\
&= K(a_{s-2},\ldots,a_t).
\end{align*}
Now, inductive application of (\ref{recursion}) to the continuants $K(b_{s-j},\ldots,b_t)$, $3 \leq j \leq s-1$, yields $K(\mathbf{b}) \geq K(\mathbf{a})$, since $a_{s-j} = b_{s-j}$ in this range.

By repeating cases (i) and (ii) as appropriate, we will eventually reach a sequence of partial quotients bounded by $B$, and at each stage we never decrease the corresponding continuant.  The result therefore follows.
\end{proof}

It remains to find a bound on $K(B,\ldots,B)$.

\begin{lemma} If $B \geq 1$, $K(\underbrace{B,\ldots,B}_m) \leq \left (\frac{1}{2} (B + \sqrt{B^2+4})\right)^{m+1}$.
\end{lemma}
\begin{proof} We proceed by induction.  The case $m=0$ is trivial.  Suppose it is true for all $m<M$.  Then, by (\ref{recursion}),
\begin{align*}
K(\underbrace{B,\ldots,B}_M) &= B K(\underbrace{B,\ldots,B}_{M-1}) + K(\underbrace{B,\ldots,B}_{M-2}) \\
& \leq B \left (\frac{1}{2} (B + \sqrt{B^2+4})\right)^M + \left (\frac{1}{2} (B + \sqrt{B^2+4})\right)^{M-1} \\
& \leq \left (\frac{1}{2} (B + \sqrt{B^2+4})\right)^{M-1} \left (\frac{1}{2} B^2 + \frac{1}{2}B \sqrt{B^2+4} + 1 \right) \\
& = \left (\frac{1}{2} (B + \sqrt{B^2+4})\right)^{M+1}.
\end{align*}
\end{proof}

\vskip 30pt

\section{The Density of $\mathcal{F}(2)$} \label{F2}

\begin{cor} There is a constant $C$ and a subset $S$ of the positive integers such that $\log |S \cap [n]| / \log n \geq \log 2/\log(1+\sqrt{2}) - o(1) \approx 0.786$ so that, for each $n \in S$, there exists a $k \in [n]$, $(k,n)=1$ so that $k/n$ has partial quotients bounded in average by $2$.
\end{cor}
\begin{proof} Let $U$ be the set of all reduced fractions $p/q$, $1 < p < q$, whose partial quotients $\mathbf{a} = (a_1,a_2,\ldots,a_m)$ are bounded in average by $2$ and such that $\mathbf{a}^\prime = (a_2,\ldots,a_m)$ is bounded in average by $2$.  The number of such $\mathbf{a}$ with $K(\mathbf{a}) \leq n$ is at least twice the number of sequences $\mathbf{a}^\prime = (a_2,\ldots,a_m)$ bounded in average by $2$ with $K(\mathbf{a}^\prime) \leq n/3$, because, if $[\mathbf{a}^\prime] = p/q$, then $K(\mathbf{a}) = a_1q+p \leq 3 K(\mathbf{a}^\prime) \leq n$.  (The fact that $\mathbf{a}^\prime$ is bounded in average by $2$ implies that $[1,\mathbf{a}^\prime]$ and $[2,\mathbf{a}^\prime]$ are also.)  Then, since every rational has at most two representations as a continued fraction, the number of elements of $U$ whose denominator is $\leq n$ is at least $\bar{S}_{n/3}(2)$, which is at least $n^{2\log 2/\log(1+\sqrt{2})-o(1)}$.  Let $S$ be the set of denominators of fractions appearing in $U$.  If $p/q=[\mathbf{a}]$ is in $U$, then $[\mathbf{a}^\prime] = (q-a_1p)/p$, so $p$ is the continuant of a sequence whose partial quotients are bounded in average by $2$.  Therefore, $\bar{S}_{n/3}(2) \leq |S \cap [n]|^2$, and we may conclude that $\log |S \cap [n]|/\log n \geq \log 2/\log(1+\sqrt{2}) - o(1)$.
\end{proof}

Attempts by the author to find a generalization of the above result to $\mathcal{F}(B)$ by applying much more careful counting arguments when $B>2$ have failed thus far.  It would also be interesting to (i), calculate the Hausdorff dimension of the set of reals in $[0,1)$ whose partial quotients are bounded in average by $B$, and (ii), draw a connection, similar to that of the ``uniform'' case, between this quantity and the asymptotic density of $\mathcal{F}(B)$.

AMS Classification Numbers: 11K50; 11K38

\end{document}